\documentclass[a4paper,10pt]{article}
\usepackage[utf8]{inputenc}
\usepackage{color}

\newcommand{\bfr}{\mathbf{r}}

\usepackage[dvips]{graphicx}
\usepackage{comment}

\newcommand{\nV}{\Psi}
\newcommand{\Li}{\mathrm{Li}}
\newcommand{\kbt}{K_B T_L}

\newcommand{\hw}{\hbar \omega_{p}}
\newcommand{\mass}{m^{*}}

\newfont{\iams}{msbm9}

\newcommand{\ot}{\frac{1}{2}}

\newcommand{\bEd}{\underline{\mathsf{E}}}
\newcommand{\cd}[1]{c_{#1}}

\title{Stochastic Galerkin Methods for the Boltzmann-Poisson system}
\author{Jos\'e A. Morales Escalante$^+$, Clemens Heitzinger$^*$}

\date{$^+$The University of Texas at San Antonio -- Departments of Mathematics and Physics \& Astronomy\\
  $^*$TU Wien -- Department of Mathematics and Geoinformation}

\begin{document}

\maketitle

\begin{abstract}
We study uncertainty quantification for a Boltzmann-Poisson system that models electron transport in semiconductors 
and the physical collision mechanisms over the charges.
We use the stochastic Galerkin method in order to handle the randomness associated with the problem. 
The main uncertainty in the Boltzmann equation concerns the initial conditions for a large number of particles, 
which is why the problem is formulated in terms of a probability density in phase space. 
The second source of uncertainty, directly related to the quantum nature of the problem, is the collision operator, 
as its structure in this semiclassical model comes from the quantum scattering matrices operating on the wave function associated 
to the electron probability density. Additional sources of uncertainty are transport, boundary data, etc. 
In this study we choose first the phonon energy as a random variable, since its value influences the energy jump appearing in the collision integral for electron-phonon scattering. 
Then we choose the lattice temperature
as a random variable, since it defines the value of the collision
operator terms in the case of electron-phonon scattering by being a
parameter of the phonon distribution. The random variable for this
case is a scalar then. Finally, we present our numerical simulations. We calculate with our stochastic Discontinuous Galerkin methods the uncertainty in kinetic moments such as density, mean energy, current, etc.\ associated to a possible physical temperature variation (assumed to follow a uniform distribution) in the lattice environment, as this uncertainty in the temperature is propagated into the electron PDF. Our results then let us predict in a real world problem setting the impact that possible variations in the lab conditions or limitations in the mathematical model (such as assumption of a constant phonon energy) will have over the calculated uncertainty in the behavior of  electronic devices. 
\end{abstract}

\section{Introduction}

\subsection{The Deterministic Boltzmann-Poisson System}

Electronic transport in semiconductors is a problem that, although definitely quantum mechanical in nature, can be approximated up to a certain point
by semiclassical models featuring quantum corrections. In the semiclassical modelling, even with deterministic laws of motion, the number of electric charge carriers
$N \gg 1$ is of the order of the Avogadro number. The consequence is that a statistical model on the semiclassical scale is extremely 
adequate due to the large number of particles (i.e., the charge
carriers, which are electrons in our problem), because it is virtually impossible to know
exactly the initial conditions of positions and momentums for all the particles (on top of quantum considerations such as the Uncertainty Principle which indicate that knowing these initial conditions exactly is completely impossible). Therefore, uncertainty in the initial condition is naturally linked to the
essence
of the electron transport problem
due to its many-carriers nature, even 
under a semiclassical approximation of 
this quantum problem, which consequently requires a statistical formulation, provided then by a particle density mechanics approach in terms of a 
probability density function in phase space. This probabilistic
formulation is given precisely by the Boltzmann-Poisson (BP) semiclassical model for collisional electronic transport.

The Boltzmann-Poisson (BP) system for electron transport on a single conduction energy band has the form
\begin{equation}
\frac{\partial f}{\partial t} + \frac{1}{\hbar} \nabla_{\vec{k}} \,
\varepsilon(\vec{k}) \cdot \nabla_{\vec{x}} f -
 \frac{q}{\hbar} \vec{E}(\vec{x},t) \cdot \nabla_{\vec{k}} f = Q(f),
\label{BE1band}
\end{equation}
\begin{equation}
\nabla_{\vec{x}} \cdot \left( \epsilon \, \nabla_{\vec{x}} V \right)
 = q \left[ \rho(\vec{x},t) - N(\vec{x}) \right] , 
\quad \vec{E} = - \nabla_{\vec{x}} V,
  \label{pois1band}
\end{equation} 
with the quantum mechanical electron group velocity 
$\frac{1}{\hbar} \nabla_{\vec{k}} \, \varepsilon (\vec{k}) $ and the 
electron density $ \rho(\vec{x},t) = \int_{\Omega_{\vec{k}}} f(\vec{x},\vec{k}, t) \, d\vec{k} $. 
The collision integral operator $Q(f)$  describes the 
scattering over the electrons,
where several mechanisms of quantum mechanical nature can be taken into account. In its full form, it enforces the Pauli Exclusion Principle by being given as
\begin{equation}
Q(f)(t, \vec{x}, \vec{k}) = \int_{\Omega_{\vec{k}'}}
\left[   S(\vec{k}' \rightarrow \vec{k}) f'(1-f)  - 
S(\vec{k} \rightarrow \vec{k}') f(1-f') \right] d \vec{k}'  .
\end{equation}
The collision scattering term $S(\vec{k}\rightarrow\vec{k}';\varepsilon(\vec{k})\rightarrow\varepsilon(\vec{k}') ) $ 
acts over $f$ in our semiclassical model  
as a scattering matrix does in a quantum description over the wave
function for a $\vec{k}$-state, i.e.,
representing the transition from a momentum $\vec{k}$ to another state $\vec{k}'$, satisfying momentum and energy conservation principles. There's then an analogy 
{\color{blue}{ $\left\langle \Psi(\vec{k}) |S| \Psi (\vec{k}')\right\rangle \leftrightarrow 
\int_{\Omega_{\vec{k}'}} S(\vec{k}'\rightarrow\vec{k}) f'(1-f) d \vec{k}'.  $ }  } 
It is important to mention that the specific form of $S(\vec{k}\rightarrow\vec{k}') $ can be derived from first-order time-dependent perturbation theory for the Schrödinger equation, by considering the perturbative Hamiltonian representing the scattering mechanisms under consideration. 
In the low density regime, however, we can relax the enforcing of the Pauli principle. In that case, the collisional 
integral operator can be approximated as being linear in $f$ and therefore having the form
\begin{equation}
Q(f) = \int_{\Omega_{\vec{k}}} \left[ S(\vec{k}', \vec{k}) f(t, \vec{x}, \vec{k}') - S(\vec{k}, \vec{k}') f(t, \vec{x}, \vec{k}) \right] d \vec{k}',
\label{ope_coll}
\end{equation}
where $S(\vec{k},\vec{k}')$ is the scattering kernel
representing non-local interactions of electrons with a background
density distribution.
For example, in the case of silicon,  one of the most important collision mechanisms are electron-phonon scatterings due to lattice vibrations of the crystal, which are modeled by  acoustic (assumed elastic) and optical (non-elastic) non-polar modes, the latter with a single frequency $\omega_{p}$ (assumed constant), as in
\begin{eqnarray}
S(\vec{k}, \vec{k}') & = & (n_{q} + 1) \, K \, \delta(\varepsilon(\vec{k}') - \varepsilon(\vec{k}) + \hw) \nonumber
\\
&&
\mbox{} + n_{q} \, K \, \delta(\varepsilon(\vec{k}') - \varepsilon(\vec{k}) - \hw) +  K_{0} \, \delta(\varepsilon(\vec{k}') - \varepsilon(\vec{k})),
\label{Skkarrow}
\end{eqnarray}
where $K$ and $K_{0}$ are material constants for silicon. 
The symbol $\delta$ indicates the usual Dirac delta distribution, derived under the
well-known Fermi's Golden Rule approximation
in time-dependent perturbation theory
\cite{lundstrom_2000}.
The constant $n_q$ is related to the phonon occupation factor
\begin{equation}
n_q(\hw) = \left[ \exp \left( \frac{\hw}{K_B T_L} \right) - 1 \right] ^{-1}, 
\end{equation}
where $K_B$ is the Boltzmann constant and $T_L = 300 K$ is the lattice temperature.

\subsection{Main Uncertainties of the Boltzmann-Poisson Model}

We can summarize the main uncertainties of the Boltzmann-Poisson model for electron transport in semiconductors as follows.
\begin{enumerate}
\item The initial conditions for and the large number of particles of
  the system, leading to a probabilistic formulation of the problem in
  terms of $f(\vec{x},\vec{k},t)$.
\item Quantum mechanical features in the Boltzmann equation, particulary in the collision operator $Q(f)$,
  based on a probabilistic description of the electron as a wavefunction~$\Psi$:
$\, \left\langle \Psi_{\vec{k}}|S|\Psi_{\vec{k}'} \right\rangle \leftrightarrow \int_{\Omega_{\vec{k}'}} S(\vec{k}'\rightarrow\vec{k}) f'(1-f) d \vec{k}' $.
\item Uncertainty in the exact functional form of the energy band $\varepsilon(\vec{k})$.  This function defines both the quantum terms of transport 
$\nabla_{\vec{k}} \varepsilon(\vec{k})$ and of electron-phonon scattering in silicon
$\delta(\varepsilon(\vec{k}) - \varepsilon(\vec{k}') + l \hw )$,
 $l \in \{-1,0,+1\}$, appearing in the Boltzmann equation.
\item The lattice temperature~$T$ may fluctuate, since it is related to the environmental temperature.
\item The phonon energy $\hw$ is often  assumed to be constant, but it is known from experiments that this is not the case in general.
\item Parameters in the Poisson equation such as doping and
  permittivity may be uncertain.  These are experimental parameters,
  and since they are given by measurements they are associated with
  measurement errors.
\item The boundary conditions that connect the domain to a stochastic environment. For example, reflection at physical boundaries might not be perfectly specular but rather have a diffusive component due to roughness in these boundaries.
\end{enumerate}
{\color{red}{In summary, uncertainty in the Boltzmann-Poisson system
    is crucial due to the inherent probabilistic (many particles) and
    quantum mechanical nature of the problem itself.}}

\subsection{The Stochastic Boltzmann-Poisson System}

Because of these reasons, we consider  a stochastic Boltzmann-Poisson system of the form
\begin{equation}
\frac{\partial f}{\partial t} + \frac{1}{\hbar} \nabla_{\vec{k}} \,
\varepsilon(\vec{k},\vec{z}) \cdot \nabla_{\vec{x}} f -
 \frac{q}{\hbar} \vec{E}(\vec{x},t,\vec{z}) \cdot \nabla_{\vec{k}} f = Q(f)(t,\vec{x},\vec{k},\vec{z}),
\label{SBE1band}
\end{equation}
\begin{equation}
\nabla_{\vec{x}} \cdot \left( \epsilon \, \nabla_{\vec{x}} V \right)
 = q \left[ \rho(\vec{x},t) - N(\vec{x}) \right] , 
\quad \vec{E} = - \nabla_{\vec{x}} V,
  \label{Spois1band}
\end{equation} 
with $f(t,\vec{x},\vec{k},\vec{z})$ being the probability density
function now depending on an additional random parameter vector
$\vec{z}$.  The components of our random vector will be associated
with the sources of uncertainty abovementioned, and they will be
explicitly restated when we describe the stochastic Galerkin method
for the Boltzmann-Poisson system.

\subsection{Previous Work on Uncertainty Quantification for Boltzmann
  Models via Stochastic Galerkin Methods}

In addition to the classical references for the stochastic Galerkin
(SG) method such as Wiener's polynomial chaos \cite{10.2307/2371268}, 
Ghanem and Spanos \cite{Ghanem:1991:SFE:103013}, 
Xiu and Karniadakis \cite{doi:10.1137/S1064827501387826}
etc., 
SG for Boltzmann equations in particular  has recently been developed mainly by Shi Jin and his collaborators. It is a very active research area in the kinetic-theory community. 
The first paper that considered the use of SG for the Boltzmann equation in the context of gases was written by 
Hu and Jin \cite{HU2016150}.
Later on, SG was studied in the context of kinetic equations with random inputs, considering different
models such as
random linear and nonlinear Boltzmann equations, linear transport equations, and Vlasov-Poisson-Fokker-Planck equations. An overall view of the advances in the discipline for these equations can be found in the review paper \cite{Hu2017},
part as well of the review book \cite{JinPareschiUQkinetic}.

More specifically, regarding SG methods
for the semiconductor Boltzmann equation, the first work related to this topic was performed by Jin and Liu
\cite{doi:10.1137/15M1053463}. 
They consider in their model a collision operator whose scattering kernel term is bounded above and below. 
Although uncertainties can possibly come from the collision operator,
the electric potential, initial data, or boundary data, the collision
operator in this previous study of the semiconductor Boltzmann equation did
not consider the more physically realistic case of Dirac delta
distributions obtained by Fermi's Golden Rule for energy transitions, as
the scattering kernel was assumed to be bounded.  Furthermore, a
possible uncertainty in the electron velocity was not considered by
making the assumption of a deterministic velocity given by the
parabolic energy band model. Most importantly, the numerical study in this work
considered a random relaxation Maxwellian collision kernel, which does
not involve energy transitions in its scattering model, in addition to random initial data in
the electron density, random boundary data, a random Debye length, and
random doping parameters in the Poisson equation.  Therefore, the
randomness of the energy band in the transport and collision terms and
the uncertainty related to a collision operator that uses Dirac delta
distributions due to Fermi's Golden Rule, as it is the case for
electron-phonon scattering, remain as crucial topics yet to be studied
for the understanding of uncertainty quantification in collisional
electron transport in semiconductors via stochastic methods.

Our methodology will then be to study variables related to the
electron-phonon collision operator as random in the SG method for the
Boltzmann-Poisson model of electrons in semiconductors.
We first choose as one of those variables the lattice temperature, as
it is involved in the phonon distribution as a parameter.  The
dimensional cost is minimal as the temperature is a scalar.

We describe the structure of the rest of this paper as follows.  In
Section~2 we discuss how the SG method handles the uncertainties
arising in the Boltzmann-Poisson system. Then we consider uncertainty
quantification by SG first for the phonon energy being a random
variable and then the lattice temperature being a random
variable too. Section~3 describes the numerics of the deterministic
Boltzmann-Poisson system solved by discontinuous-Galerkin (DG)
methods. Then Section~4 covers in more detail the stochastic
discontinuous Galerkin (SDG) method for the Boltzmann-Poisson system for the
case of a random lattice temperature. In Section~5, the conclusions
are drawn.

\section{Stochastic Galerkin Method for the Boltzmann-Poisson System}

\subsection{Description of Uncertainties}

The SG method handles the uncertainties in the Boltzmann-Poisson
system by introducing random variables $z_i$, $i \in \{1, \ldots,
7\}$, associated with the uncertainties as indicated below.

\begin{enumerate}
\item Regarding initial conditions and the large number of particles,
  the probabilistic formulation {\color{blue}{
      $f(\vec{x},\vec{k},t,z_1)$}} with random initial
  conditions {\color{blue}{
      $f(\vec{x},\vec{k},0,z_1)$}} is used.
\item Regarding quantum phenomena in the collision operator $Q(f)$,
  the probabilistic nature of the electron wavefunction $\Psi$ is
  mimicked (relaxing Pauli principle) by $\left\langle \Psi_{\vec{k}'}|S|\Psi_{\vec{k}} \right\rangle
  \leftrightarrow \int_{\Omega_{\vec{k}}}
  {\color{blue}{ S(\vec{k} \rightarrow
      \vec{k}',z_2) f(\vec{x},\vec{k},t,z_2) }} d\vec{k}'$.
\item The uncertainty in the energy band structure is described as
  {\color{blue}{ $\varepsilon(\vec{k},z_3)$}}.  This function
  defines both quantum terms of electron velocity {\color{blue}{
      $\nabla_{\vec{k}} \varepsilon(\vec{k},z_3)/\hbar$}} and of
  electron-phonon scattering in silicon {\color{blue}{
      $\delta(\varepsilon(\vec{k},z_3) -
      \varepsilon(\vec{k}',z_3) + l (\hw + z_4) ),$ $l \in \{-1,0,+1\}$.}}
\item The lattice temperature is written as {\color{blue}{ $T + z_4
      $}} as it may change due to fluctuations in the environment, but
  it is assumed to be constant in the model.
\item The phonon energy {\color{blue}{ $\hw + z_5 $}} is approximated as
  constant in the model, but experiments show that it is nonconstant in general. 
\item In the Poisson equation, the doping concentration is written as
  {\color{blue}{ $N_D + z_6 $}} and the permittivity as {\color{blue}{
      $\epsilon + z_6 $}}.
\item The boundary conditions are described as {\color{blue}{
      $f_B(\vec{x},\vec{k},t,z_7)|_{\partial\Omega}$}}, etc.
\end{enumerate}


\subsection{Stochastic Galerkin Method for the Boltzmann-Poisson
  System with the Phonon Energy as a Random Variable}

In this section, we assume that the only uncertainty stems from the
phonon energy model.  Randomness in the phonon energy is physically relevant because it is known to be strictly speaking non-constant, and it is also a good
starting point as a scalar random variable. 
Therefore, we replace $\hw$ in the deterministic equation by $\hw +
z$.  Then the phonon occupation as a function of the energy
becomes
\begin{equation}
n_q(\hw,z) = \left[ \exp \left( \frac{\hw + z}{K_B T_L} \right) - 1 \right] ^{-1}.
\end{equation}
This introduces randomness in the collision operator as well, leading to
\begin{eqnarray*}
S(\vec{k}, \vec{k}',z) & = & \left[n_{q}(\hw,z) + 1\right] \, K \, \delta(\varepsilon(\vec{k}') - \varepsilon(\vec{k}) + \hw + z) 
\\
+  K_{0} \, \delta(\varepsilon(\vec{k}') - \varepsilon(\vec{k}) ) 
&+&
  n_{q}(\hw,z) \, K \, \delta(\varepsilon(\vec{k}') - \varepsilon(\vec{k}) - \hw - z) .
\end{eqnarray*}
We consider two cases.  First we devise a stochastic Galerkin
algorithm using a distributional-derivative approximation with respect
to the random variable. Then we consider the fully general case of the
random variable in the collision operator without any
distributional-derivative approximation in the random space. The
algorithms are described in detail in the following.

\subsubsection{Random Phonon Energy Using a Distributional Derivative Approximation}

We consider the case $n :=1$, $P :=1$, $K:=\dim(\mathcal{P}^n_P):= 2$,
and approximate the density as (with
$\vec{p}=\hbar\vec{k}$ being the crystal momentum)
\begin{equation}
  f(t,\vec{x},\vec{p},z) \approx \sum_{k=1}^{2} \alpha_k(t,\vec{x},\vec{p}) \Psi_k(z)  = \alpha(t,\vec{x},\vec{p})\cdot \Psi(z) = (\alpha_1,\alpha_2) \cdot (\Psi_1,\Psi_2).
\end{equation}
Then the Boltzmann equation reads
\begin{equation}
\partial_t \alpha + \vec{v} \cdot \nabla_{\vec{x}} \alpha + F \cdot \nabla_{\vec{p}} \alpha = Q(\alpha),
\end{equation}
\begin{eqnarray}
Q(\alpha) &=&  \int_{\Omega_{\vec{p}}} B(\vec{p},\vec{p'}) [M(\vec{p})\alpha(\vec{p'}) - M(\vec{p'})\alpha(\vec{p})] d\vec{p'},\\
B_{ij}(\vec{p},\vec{p'}) &=&  \int_{I_z} \sigma(\vec{p},\vec{p'},z) \Psi_i(z) \Psi_j(z) \pi(z) dz\\
  &=& \sigma_0(\vec{p},\vec{p'}) \delta_{ij} +  \int_{I_z} \partial_z{\sigma}(\vec{p},\vec{p'},z)|_{z=0}z \Psi_i(z) \Psi_j(z) \pi(z) dz,\\
\sigma(\vec{p},\vec{p'},z) &=& \sigma_0(\vec{p},\vec{p'}) + \tilde{\sigma}_1(\vec{p},\vec{p'})z = \sigma|_{z=0} + \partial_z{\sigma}(\vec{p},\vec{p'},z)|_{z=0}z.
\end{eqnarray}
The scattering cross section is then written as
\begin{equation}
\sigma(\vec{p},\vec{p'},z) = \sigma_0(\vec{p},\vec{p'}) + \tilde{\sigma}_1(\vec{p},\vec{p'})z ,
\end{equation}
where
\begin{equation}
\sigma_0  = M^{-1}[K_0 \delta(\varepsilon - \varepsilon') +(e^{\hw} -1)^{-1}K(e^{\hw}\delta(\varepsilon - \varepsilon' + \hw)
+ \delta(\varepsilon - \varepsilon' - \hw) ) ].
\end{equation}
Here
$\tilde{\sigma}_1(\vec{p},\vec{p'}) = \partial_z{\sigma}(\vec{p},\vec{p'},z)|_{z=0} $ is a
distributional derivative with respect to $z$, and we find
\begin{eqnarray*}
\tilde{\sigma}_1(\vec{p},\vec{p'})  &=& M^{-1}K \partial_z \left\lbrace [ 1 + (e^{\hw+z}-1)^{-1}] \delta(\varepsilon - \varepsilon' + \hw + z) \right. \\
&& \left.
{} + (e^{\hw+z} -1)^{-1} \delta(\varepsilon - \varepsilon' - \hw -z ) \right\rbrace |_{z=0}.
\end{eqnarray*}
Using $\delta'[\phi] = - \delta[\phi']$ and phonon distribution properties, we have
\begin{eqnarray*}
\tilde{\sigma}_1(\vec{p},\vec{p'})  &=& M^{-1}K  \left\lbrace [ 1 + (e^{\hw}-1)^{-1}] \partial_z|_{0}\delta(\varepsilon - \varepsilon' + \hw + z) \right. \\
&& 
{} + (e^{\hw} -1)^{-1} \partial_z|_{0}\delta(\varepsilon - \varepsilon' - \hw -z ) \\
&& \left.
{} -\frac{e^{\hw}}{(e^{\hw} -1)^2} [\delta(\varepsilon - \varepsilon' + \hw ) + \delta(\varepsilon - \varepsilon' - \hw )] 
\right\rbrace 
\end{eqnarray*}
and
\begin{eqnarray*}
B_{ij}(\vec{p},\vec{p'})
&=& \sigma_0(\vec{p},\vec{p'}) \delta_{ij}  + M^{-1}K  \int_{I_z} dz \Psi_i(z) \Psi_j(z) \pi(z) z
 \cdot \\ 
&& \cdot
\left\lbrace [ 1 + (e^{\hw}-1)^{-1}] \partial_z|_{0}\delta(\varepsilon - \varepsilon' + \hw + z) \right. \\
&& 
{} + (e^{\hw} -1)^{-1} \partial_z|_{0}\delta(\varepsilon - \varepsilon' - \hw -z ) \\
&& \left.
{} -\frac{e^{\hw}}{(e^{\hw} -1)^2}[\delta(\varepsilon - \varepsilon' + \hw ) + \delta(\varepsilon - \varepsilon' - \hw )] 
\right\rbrace   \\
&=& \sigma_0(\vec{p},\vec{p'}) \delta_{ij}  - M^{-1}K  
 \cdot \\ 
&& 
\cdot \left\lbrace [ 1 + (e^{\hw}-1)^{-1}] \partial_z[\Psi_i(z) \Psi_j(z) \pi(z) z]\chi|_{z= - (\varepsilon - \varepsilon' + \hw)} \right. \\
&& 
{}+ (e^{\hw} -1)^{-1} \partial_z[\Psi_i(z) \Psi_j(z) \pi(z) z]\chi|_{z= + (\varepsilon - \varepsilon' - \hw)} \\
&& \left.
{}+\frac{ e^{\hw} \int_{I_z} dz \Psi_i(z) \Psi_j(z) \pi(z) z }{ (e^{\hw} -1)^2} \sum_{\pm} \delta(\varepsilon - \varepsilon' \pm \hw )   
\right\rbrace.
\end{eqnarray*}

If we assume the example $\pi(z) := \frac{e^{-\frac{z^2}{2}}}{\sqrt{2}}$,
$\Psi_1 := 1$, and $\Psi_2 := 2z$, we obtain
\begin{eqnarray*}
B(\vec{p},\vec{p'})
&=& \sigma_0(\vec{p},\vec{p'}) \delta_{ij}  - M^{-1}K  
 \times \\ 
&& 
\left\lbrace \left[ 1 + \frac{1}{e^{\hw}-1}\right] \partial_z\left[\left(
\begin{array}{*{2}c}
 1 & 2z \\
 2z & 4z^2
\end{array}
\right)
 \pi(z) z\right]\chi|_{z= - (\varepsilon - \varepsilon' + \hw)} \right. \\
&& 
{} + (e^{\hw} -1)^{-1} \partial_z\left[\left(
\begin{array}{*{2}c}
 1 & 2z \\
 2z & 4z^2
\end{array}
\right) \pi(z) z\right]\chi|_{z= + (\varepsilon - \varepsilon' - \hw)} \\
&& \left.
{} +\frac{ e^{\hw} \sum_{\pm} \delta(\varepsilon - \varepsilon' \pm \hw )  }{ (e^{\hw} -1)^2}
\left(
\begin{array}{*{2}c}
0 & 1/2 \\
1/2 & 0
\end{array}
\right)
\right\rbrace   \\
&=& \sigma_0(\vec{p},\vec{p'}) \mathcal{I}  - M^{-1}K  
 \times \\ 
&& 
\left\lbrace    \frac{e^{\hw}}{e^{\hw}-1}  \left(
\begin{array}{*{2}c}
 1-z^2 & 2z(2-z^2) \\
 2z(2-z^2) & 2z^2(3-z^2)
\end{array}
\right)
 \pi \chi|_{z= - (\varepsilon - \varepsilon' + \hw)} \right. \\
&& 
{} + \frac{1}{e^{\hw}-1}     \left(
\begin{array}{*{2}c}
 1-z^2 & 2z(2-z^2) \\
 2z(2-z^2) & 2z^2(3-z^2)
\end{array}
\right)  \pi \chi|_{z= + (\varepsilon - \varepsilon' - \hw)} \\
&& \left.
{} +\frac{ e^{\hw}  }{ (e^{\hw} -1)^2} \sum_{\pm} \delta(\varepsilon - \varepsilon' \pm \hw )
\left(
\begin{array}{*{2}c}
0 & 1/2 \\
1/2 & 0
\end{array}
\right)
\right\rbrace,\\
Q(\alpha) &=&  \int_{\Omega_{\vec{p}}} B(\vec{p},\vec{p'}) \left[M(\vec{p})\alpha(\vec{p'}) - M(\vec{p'})\alpha(\vec{p})\right] d\vec{p'}.
\end{eqnarray*}

One could use as another example a sharper Gaussian with support mostly concentrated around the central value of zero fluctation, in order to give such a low probability to energy values $\hbar\omega_p + z< 0$ that virtually the probability density of having $z< -\hbar\omega_p$ would be numerically zero in a computational implementation.

\subsubsection{Random Phonon Energy Without Approximation}

Next, we consider again a random phonon energy, but now use the
collision scattering term without approximation by distributional
derivatives in the random space.
We consider the case
$n :=1$, $P :=1$, $K:=\dim(\mathcal{P}^n_P) := 2$,
and approximate the density again as
\begin{equation}
f(t,\vec{x},\vec{p},z) \approx \sum_{k=1}^{2} \alpha_k(t,\vec{x},\vec{p}) \Psi_k(z)  = \alpha(t,\vec{x},\vec{p})\cdot \Psi(z) = (\alpha_1,\alpha_2) \cdot (\Psi_1,\Psi_2).
\end{equation}
Then the Boltzmann equation reads
\begin{equation}
\partial_t \alpha + v \cdot \nabla_{\vec{x}} \alpha + F \cdot \nabla_{\vec{p}} \alpha = Q(\alpha),
\end{equation}
\begin{eqnarray}
Q(\alpha) &=&  \int_{\Omega_p} B(\vec{p},\vec{p'}) [M(\vec{p})\alpha(\vec{p'}) - M(\vec{p'})\alpha(\vec{p})] d\vec{p'}, \\
B_{ij}(\vec{p},\vec{p'}) &=&  \int_{I_z} \sigma(\vec{p},\vec{p'},z) \Psi_i(z) \Psi_j(z) \pi(z) dz,
\end{eqnarray}
\begin{equation}
\sigma(\vec{p},\vec{p'},z) = 
\frac{
K_0 \delta(\varepsilon - \varepsilon') +
K \frac{
e^{\beta(\hw + z)}\delta(\varepsilon - \varepsilon' + \hw + z)
+ \delta(\varepsilon - \varepsilon' - \hw - z)  }{e^{\beta(\hw + z)} -1}
}{M(\vec{p})}
\end{equation}
with $\beta=(\kbt)^{-1}$. We then have
\begin{eqnarray*}
B_{ij}(\vec{p},\vec{p'})
&=& \frac{ K_0\delta(\varepsilon-\varepsilon')\delta_{ij} +  K  \int_{I_z} dz \Psi_i(z) \Psi_j(z) \pi(z) 
\left\lbrace 
 \frac{
e^{\beta(\hw + z)}\delta(\varepsilon - \varepsilon' + \hw + z)
+ \delta(\varepsilon - \varepsilon' - \hw - z)  }{e^{\beta(\hw + z)} -1}
\right\rbrace  }{M(\vec{p})}.
\end{eqnarray*}
Therefore, we find
\begin{eqnarray*}
B_{ij}(\vec{p},\vec{p'})
&=& \frac{ K_0\delta(\varepsilon-\varepsilon')\delta_{ij} +  K \left( \left.  
 \frac{ \chi(z) \Psi_i(z) \Psi_j(z) \pi(z)
e^{\beta(\hw + z)} }{e^{\beta(\hw + z)} -1}
\right|_{z = \varepsilon' - \varepsilon - \hw }
+  \left.  
 \frac{\chi(z) \Psi_i(z) \Psi_j(z) \pi(z)}{e^{\beta(\hw + z)} -1}
\right|_{z = \varepsilon - \varepsilon' - \hw }
\right)  }{M(\vec{p})}
\end{eqnarray*}
with $\chi(z)$ being the characteristic function.  Furthermore, we
have
\begin{eqnarray*}
B &=& M^{-1}(\vec{p}) K_0\delta(\varepsilon-\varepsilon')I \\ 
&+& \frac{   \left.  
 \frac{ K \chi(z)  \pi(z) }{1 - e^{-\beta(\hw + z)} }
\left(
\begin{array}{*{2}c}
\Psi_1^2(z) & \Psi_1 \Psi_2 \\
\Psi_1 \Psi_2 & \Psi_2^2(z)
\end{array}
\right)
\right|_{\varepsilon' - \varepsilon - \hw }
+  \left.  
 \frac{K \chi(z) \pi(z)}{e^{\beta(\hw + z)} -1}
 \left(
\begin{array}{*{2}c}
\Psi_1^2(z) & \Psi_1 \Psi_2 \\
\Psi_1 \Psi_2 & \Psi_2^2(z)
\end{array}
\right)
\right|_{\varepsilon - \varepsilon' - \hw }
  }{M(\vec{p})} .
\end{eqnarray*}
If we assume a uniform distribution $\pi(z) = {N}/{2\beta}$ for $z \in
[-\beta/N, \beta/N]$ with $N>1$, or equivalently $\pi(w) = {1}/{2}$ by
the scaling $w = Nz/\beta$ for $w \in [-1, 1]$, with the associated
Legendre polynomials $\Psi_1 = 1$ and $\Psi_2(w) = w$, we obtain
\begin{eqnarray*}
B &=& M^{-1}(\vec{p}) K_0\delta(\varepsilon-\varepsilon')I \\ 
&+& \frac{   \left.  
 \frac{ K \chi(z) /2 }{1 - e^{-\beta(\hw + z)} }
\left(
\begin{array}{*{2}c}
1 & \frac{ Nz}{\beta} \\
\frac{ Nz}{\beta} & (\frac{ Nz}{\beta})^2
\end{array}
\right)
\right|_{z = \varepsilon' - \varepsilon - \hw }
{}+  \left.  
 \frac{K \chi(z) /2 }{e^{\beta(\hw + z)} -1}
 \left(
\begin{array}{*{2}c}
1 & \frac{ Nz}{\beta} \\
\frac{ Nz}{\beta} & (\frac{ Nz}{\beta})^2
\end{array}
\right)
\right|_{z = \varepsilon - \varepsilon' - \hw }
  }{M(\vec{p})},
\end{eqnarray*}
\begin{eqnarray*}
Q(\alpha) &=&  \int_{\Omega_{\vec{p}}} B(\vec{p},\vec{p'}) \left[M(\vec{p})\alpha(\vec{p'}) - M(\vec{p'})\alpha(\vec{p})\right] d\vec{p'}.
\end{eqnarray*}


\subsection{Stochastic Galerkin Method for the Boltzmann-Poisson
  System with the Lattice Temperature as a Random Variable}

In this example, we assume that the only uncertainty in our problem
stems from the lattice temperature.  Randomness in the lattice
temperature is motivated by physical reasons, as the temperature in
the material or in its environment often fluctuates.  The random
variable is scalar, and randomness is introduced in the collisions,
but now outside the argument of the Dirac delta distributions associated
with Fermi's Golden Rule.


Therefore, the term $\kbt$ in the deterministic equation is replaced
by $\kbt + z^*$, or equivalently, the term $\beta := (\kbt)^{-1}$ in
the deterministic equation is replaced by $\beta + z$.  This
introduces randomness in the phonon occupation as a function the
energy, yielding
\begin{equation}
n_q(\hw,z) = \left[ \exp \left( \frac{\hw }{K_B T_L  + z^*} \right) - 1 \right] ^{-1} 
=
\left[ e^{(\beta + z)\hw} - 1 \right] ^{-1}.
\end{equation}
Additionally, randomness in the collision operator model is introduced
as well. We have
\begin{eqnarray*}
S(\vec{k}, \vec{k}',z) & = & \left[n_{q}(\hw,z) + 1\right] \, K \, \delta(\varepsilon(\vec{k}') - \varepsilon(\vec{k}) + \hw )
\\
{} +  K_{0} \, \delta(\varepsilon(\vec{k}') - \varepsilon(\vec{k}) ) 
&+&
  n_{q}(\hw,z) \, K \, \delta(\varepsilon(\vec{k}') - \varepsilon(\vec{k}) - \hw ) .
\end{eqnarray*}
Noticing that the randomness is just in the coefficients related to
the phonon density and not inside the arguments of the delta
distributions, we equivalently have
\begin{eqnarray*}
S(\vec{k}, \vec{k}',z) & = &  K \, \frac{ e^{(\beta + z)\hw} }{  e^{(\beta + z)\hw} - 1   } \,  \delta(\varepsilon(\vec{k}') - \varepsilon(\vec{k}) + \hw ) 
\\
{} +  K_{0} \, \delta(\varepsilon(\vec{k}') - \varepsilon(\vec{k}) ) 
&+&
K \, \frac{ 1 }{  e^{(\beta + z)\hw} - 1   } \,  \delta(\varepsilon(\vec{k}') - \varepsilon(\vec{k}) - \hw ) .
\end{eqnarray*}


In particular, we consider the random temperature $T_L + z$, set
$n :=1$, $P :=1$, $K:=\dim(\mathcal{P}^n_P) = 2$,
and approximate the density as
\begin{equation}
f(t,\vec{x},\vec{p},z) \approx \sum_{k=1}^{2} \alpha_k(t,\vec{x},\vec{p}) \Psi_k(z)  = \alpha(t,\vec{x},\vec{p})\cdot \Psi(z) = (\alpha_1,\alpha_2) \cdot (\Psi_1,\Psi_2).
\end{equation}
This yields the Boltzmann equation
\begin{equation}
\partial_t \alpha + v \cdot \nabla_{\vec{x}} \alpha + F \cdot \nabla_{\vec{p}} \alpha = Q(\alpha),
\end{equation}
\begin{eqnarray}
Q(\alpha) &=&  \int_{\Omega_p} B(\vec{p},\vec{p'}) [M(\vec{p})\alpha(\vec{p'}) - M(\vec{p'})\alpha(\vec{p})] d\vec{p'},\\
B_{ij}(\vec{p},\vec{p'}) &=&  \int_{I_z} \sigma(\vec{p},\vec{p'},z) \Psi_i(z) \Psi_j(z) \pi(z) dz
\end{eqnarray}
with the scattering cross section
\begin{equation}
\sigma(\vec{p},\vec{p'},z)  = \frac{1}{M(\vec{p})} \left(K_0 \delta(\varepsilon - \varepsilon') +
K \frac{
e^{(\beta+z)\hw}\delta(\varepsilon - \varepsilon' + \hw)
+ \delta(\varepsilon - \varepsilon' - \hw) 
}{e^{(\beta +z)\hw} -1}
 \right).
\end{equation}
Substituting $\sigma$ into $B_{ij}$ yields
\begin{eqnarray*}
B_{ij}(\vec{p},\vec{p'}) &=&  \int_{I_z} dz \Psi_i(z) \Psi_j(z) \pi(z) 
\frac{ K_0 \delta(\varepsilon - \varepsilon') +
K \frac{
e^{(\beta+z)\hw}\delta(\varepsilon - \varepsilon' + \hw)
+ \delta(\varepsilon - \varepsilon' - \hw) 
}{e^{(\beta +z)\hw} -1}
}{M(\vec{p})} \\
 &=&  
M^{-1}(\vec{p})K_0 \delta(\varepsilon - \varepsilon') \delta_{ij}  +
 \int_{I_z} dz \Psi_i(z) \Psi_j(z) \pi(z) 
\frac{ 
K \frac{
e^{(\beta+z)\hw}\delta(\varepsilon - \varepsilon' + \hw)
+ \delta(\varepsilon - \varepsilon' - \hw) 
}{e^{(\beta +z)\hw} -1}
}{M(\vec{p})}.
\end{eqnarray*}
This expression simplifies to
\begin{eqnarray*}
B_{ij}(\vec{p},\vec{p'}) &=&  
\frac{
K_0 \delta(\varepsilon - \varepsilon') \delta_{ij}  +
K \left( \delta(\varepsilon - \varepsilon' + \hw)
\int_{I_z}  
\frac{  dz \Psi_i \Psi_j \pi
}{1 - e^{-(\beta +z)\hw} }
+
 \delta(\varepsilon - \varepsilon' - \hw)
\int_{I_z}  
\frac{dz \Psi_i \Psi_j \pi}{e^{(\beta +z)\hw} -1}
\right)
}{M(\vec{p})}
\end{eqnarray*}
and
\begin{eqnarray*}
B_{ij}(\vec{p},\vec{p'}) &=&  
{M(\vec{p})}^{-1}
\left[
K_0 \delta(\varepsilon - \varepsilon') \delta_{ij}  +
K \left( \delta(\varepsilon - \varepsilon' + \hw) C^+_{ij}
+
 \delta(\varepsilon - \varepsilon' - \hw) C^-_{ij}
\right)
\right]
\end{eqnarray*}
with the coefficients
\begin{eqnarray*}
C^-_{ij} &=& \int_{I_z}  dz \Psi_i(z) \Psi_j(z) \pi(z)
\frac{1}{e^{(\beta +z)\hw} -1}
=
\int_{I_z}  dz \Psi_i(z) \Psi_j(z) \pi(z) \,
n_q(\hw, \beta +z) ,
\\
C^+_{ij} &=& 
\int_{I_z} dz \Psi_i(z) \Psi_j(z) \pi(z)  
\left( 1 +
\frac{1}{e^{(\beta +z)\hw} - 1 }
\right)
=
\int_{I_z} dz \Psi_i \Psi_j \pi  
\left( n_q + 1 
\right)
= \delta_{ij} + C^-_{ij} .
\end{eqnarray*}
Therefore the scattering operator becomes
\begin{eqnarray*}
Q(\alpha) &=&  \int_{\Omega_{\vec{p}}} 
\left(
K_0 \delta(\varepsilon - \varepsilon') I_{d}  +
K \sum_{\pm} \delta(\varepsilon - \varepsilon' \pm \hw) C^{\pm}
\right)
M(\vec{p})^{-1}
\left[ M(\vec{p})\alpha(\vec{p'}) - M(\vec{p'})\alpha(\vec{p})\right] d\vec{p'} 
\\
&=&  \int_{\Omega_{\vec{p}}} 
\left(
 \sum_{l=-1}^{1} K_l \delta(\varepsilon - \varepsilon' + l \hw) C^{l}
\right)
M(\vec{p})^{-1}
\left[ M(\vec{p})\alpha(\vec{p'}) - M(\vec{p'})\alpha(\vec{p})\right] d\vec{p'} 
\end{eqnarray*}
with $K_{-1} = K = K_{+1}$ and $C^0_{ij} = \delta_{ij}$ being the
identity matrix.

It is important to note that a Gaussian distribution is not
appropriate in this example, as there would arise a singularity in the
integrals when the temperature (in energy units) becomes zero.  We
hence assume a uniform distribution $\pi(z) = {N}/{2\beta}$ for $z \in
[-\beta/N, \beta/N]$ with $N>1$, or equivalently $\pi(w) = {1}/{2}$ by
the scaling $w = Nz/\beta$ for $w \in [-1, 1]$, with the associated
Legendre polynomials $\Psi_1 = 1$ and $\Psi_2(w) = w$.  This means
that $f \approx \alpha_1 + w \alpha_2 $.  Therefore we find
\begin{eqnarray*}
C^- &=& \int_{-\beta/N}^{\beta/N}  
\frac{dz \, {N}/{2\beta} }{e^{(\beta +z)\hw} -1}
\left(
\begin{array}{*{2}c}
1 & Nz/\beta \\
Nz/\beta & (Nz/\beta)^2
\end{array}
\right)
= \frac{1}{2}
\int_{-1}^{1}  
\frac{  \left(
\begin{array}{*{2}c}
1 & w \\
w & w^2
\end{array}
\right)
dw }{e^{\beta\hw(1 + w/N) } -1}.
\end{eqnarray*}
The analytic values of these integrals are
\begin{eqnarray*}
\int \frac{ dx}{\exp(A + B x) - 1}  &=& \frac{\log(1 - e^{A + B x})}{B} - x + ct ,
\\
\int \frac{ x dx}{\exp(A + B x) - 1}  &=& \frac{ \Li_2(e^{A + B x})}{B^2} + \frac{x \log(1 - e^{A + B x})}{B}  - \frac{ x^2}{2} + ct,
\\
\int \frac{ x^2 dx }{\exp(A + B x) - 1}  &=& \frac{-2 \Li_3(e^{A + B x})}{B^3} + 
\frac{2 x \Li_2(e^{A + B x})}{B^2} + \frac{x^2 \log(1 - e^{A + B x})}{B} - \frac{ x^3}{3} + ct
\end{eqnarray*}
with $A = \beta\hw$, $B = \beta\hw/N$, and $\Li_n(x)$ being the
polylogarithm functions.  Furthermore, we can evaluate these formulas
to obtain $C^-$ explicitly in the form
\begin{eqnarray*}
C^- &=& \frac{1}{2}
\left(
\begin{array}{*{2}c}
\left. \frac{\log(1 - e^{A + B x})}{B} - x \right|_{-1}^{1} & 
\left. \frac{ \Li_2(e^{A + B x})}{B^2} + \frac{x \log(1 - e^{A + B x})}{B}   \right|_{-1}^{1} \\
\left. \frac{ \Li_2(e^{A + B x})}{B^2} + \frac{x \log(1 - e^{A + B x})}{B}   \right|_{-1}^{1} & \left. \frac{-2 \Li_3(e^{A + B x})}{B^3} + 
\frac{2 x \Li_2(e^{A + B x})}{B^2} + \frac{x^2 \log(1 - e^{A + B x})}{B} - \frac{ x^3}{3} \right|_{-1}^{1}
\end{array}
\right),
\end{eqnarray*}
where we can omit the term $- \frac{ x^2}{2}$ in the off-diagonal
elements, since it will vanish when evaluating at $\pm 1$.

To determine concrete numbers for numerical simulations and for the
evaluation of the $C^-$ matrix, we recall that the Planck constant
divided by $2\pi$ is equal to $ \hbar = h/2\pi = 1.0546 \times
10^{-34} \,\mathrm{J\cdot s}$ and that the Boltzmann constant is equal
to $ K_B = 1.3805 \times 10^{-23} \,\mathrm{J/K}$. The mean lattice
temperature is assumed to be $T_L := 300 K = 26.85 ^\circ\mathrm{C}$.
Therefore, we have that $\kbt = 4.1415 \times 10^{-21} \,\mathrm{J} =
0.025\,849 \,\mathrm{eV}$, since $1 \,\mathrm{eV} = 1. 602\,18 \times
10^{-19} \,\mathrm{J}$.  Hence $\beta = (\kbt)^{-1} = 2.414\,584\,1 \times
10^{20} \,\mathrm{J}^{-1}$.

Moreover, the variation in the environment temperature might be of
$\pm 10 ^\circ\mathrm{C}$, resulting in a lattice temperature between
$16.85 ^\circ\mathrm{C} = 290 \,\mathrm{K}$ and $36.85
^\circ\mathrm{C} = 310 \,\mathrm{K}$.  In that case $\kbt \in
[4.003\,45, 4.279\,55]\times 10^{-21} \,\mathrm{J}$, $\beta + z \in
[2.336\,694\,3, 2.497\,845\,6]\times 10^{20} \,\mathrm{J}^{-1}$, and
$z \in [-0.077\,889\,8, 0.083\,261\,5]\times 10^{20}
\,\mathrm{J}^{-1}$.  Thus $z \in I_z$, which is $I_z \approx
[-0.080\,575\,65, 0.080\,575\,65]\times 10^{20} \,\mathrm{J}^{-1}$,
and therefore $\beta/N = 0.080\,575\,65 \times 10^{20}
\,\mathrm{J}^{-1}$ implies $N = \frac{ 2.414\,584\,1 \times 10^{20}
  \,\mathrm{J}^{-1}}{0.080\,575\,65 \times 10^{20} \,\mathrm{J}^{-1}}
= 29.966\,672\,313\,5$.  After rounding to $N := 30$, we have $z \in
[-\beta/N, \beta/N] $ with $\beta/N = 0.080\,486\,136\,66 \times
10^{20} \,\mathrm{J}^{-1}$.  Finally, since the phonon energy is $\hw
= 0.063 \,\mathrm{eV} =1.009\,373\,4\times 10^{-20} \,\mathrm{J}$, we
obtain the values of the adimensional numbers $\beta \hw =
2.437\,216\,962\,6 = A$ and $\beta \hw /N = 0.081\,240\,565\,42 = B$.

Hence we find the matrices
\begin{eqnarray*}
C^- &=& \frac{1}{2}
\left(
\begin{array}{*{2}c}
0.191\,825 & -0.005\,690\,12
\\
-0.005\,690\,12 & 0.064\,015\,1
\end{array}
\right)
=
\left(
\begin{array}{*{2}c}
0.095\,912\,5 & -0.002\,845\,06
\\
-0.002\,845\,06 & 0.032\,007\,55
\end{array}
\right),
\\
C^+ &=& 
\left(
\begin{array}{*{2}c}
1.095\,912\,5 & -0.002\,845\,06
\\
-0.002\,845\,06 & 1.032\,007\,55
\end{array}
\right)
,
\end{eqnarray*}
since $C^+_{ij} = C^-_{ij} + \delta_{ij}$.  With these matrices, the
scattering operator becomes
\begin{equation}
 Q(\alpha) =  \int_{\Omega_{\vec{p}}} 
\left(
K_0 \delta(\varepsilon - \varepsilon') I_d  +
K \sum_{\pm} \delta(\varepsilon - \varepsilon' \pm \hw) C^{\pm}
\right)
\frac{
M \alpha' - M'\alpha
}{M} d\vec{p'} ,
\end{equation}
which we can write in the form 
\begin{eqnarray*}
 Q(\alpha) &=&  \int_{\Omega_{\vec{p}}} 
\left(
K_0 \delta(\varepsilon - \varepsilon')   +
K \sum_{\pm} 
\delta(\varepsilon - \varepsilon' \pm \hw) \left[n_q + (1 \pm 1)/{2} \right]
\right) I_d
\frac{
M \alpha' - M'\alpha
}{M} d\vec{p'} 
\\
 &+&  \int_{\Omega_{\vec{p}}} 
K\left[
 \sum_{\pm} 
\delta(\varepsilon - \varepsilon' \pm \hw) \right]
  \left( C^{-} - n_q  I_d \right) 
\frac{
M \alpha' - M'\alpha
}{M} d\vec{p'} ,
\end{eqnarray*}
since $ C^+ = C^- + I_d$ and hence $C^- - n_q I_d=C^+ - (n_q + 1)
I_d$.

Therefore the term in the first row is the collision operator, as
originally written in the deterministic case, acting on each separate
band (by means of the identity matrix) without any recombination,
whereas the second term represents the recombination and diagonal
terms related to the uncertainty in the temperature associated solely
with inelastic integrals. Given the value of the constant $ n_q =
\left[ e^{\beta \hw} - 1 \right] ^{-1} = 0.095\,774\,842\,71$, we find
\begin{eqnarray*}
C^- - n_q I_d &=& 
\left(
\begin{array}{*{2}c}
0.000\,137\,657\,29 & -0.002\,845\,06
\\
-0.002\,845\,06 & -0.063\,767\,292\,71
\end{array}
\right)
=
C^+ - (n_q + 1) I_d.
\end{eqnarray*}

\section{Stochastic Galerkin Method for the Boltzmann-Poisson System
  Using Deterministic Discontinuous Galerkin Solvers}

The numerics of deterministic solvers for the Boltzmann-Poisson system
that use the discontinuous Galerkin (DG) algorithm have been studied
in \cite{MR2567861}, \cite{1593f467c520433c872fc676038d749a} for a
single PDF (one band) without randomness. We will use the
deterministic DG method for two bands (representing the $\alpha$
vector of coefficients) to solve the stochastic Galerkin system, which
contains a different kind of matrix integral collisional operator.

\subsection{Discontinuous Galerkin: the Boltzmann Equation in $\vec{k}$-Spherical Coordinates}

We perform a spherical transformation of the momentum coordinate
$\vec{k}$ taking the location of a (local) minimum of the conduction
energy band as the origin. This transformation is useful (in the
absence of Umklapp effects), because in low energy limits (i.e., for
small potential bias) the conduction band energy scales as the square
of the momentum norm, and hence the radial coordinate is an energy
variable. We then have
$$
\vec{k} = \frac{\sqrt{2 m^* \kbt}}{\hbar} \, \sqrt{r} \left( \mu,
\sqrt{1 - \mu^{2}} \cos \varphi,  \sqrt{1 - \mu^{2}} \sin
\varphi\right),
$$ 
$$
r \geq 0,\qquad
\mu \in [-1,1],\qquad
\varphi \in [-\pi, \pi].
$$

The variable~$r$ is proportional to the energy for small biases in the
parabolic band approximation, assuming the same effective mass in all
three Cartesian momentum directions.  Due to this momentum coordinate
tranformation, we have to weight the PDF coefficients by the Jacobian
of the $\vec{k}$-transformation, specifically for the computation of moment
integrals over the $\vec{k}$-space. We then obtain a transformed PDF in the
phase space $(\vec{x},r,\mu,\varphi)$ given by
$$\Phi(t,\vec{x},r,\mu,\varphi) = \frac{\sqrt{r}}{2} \alpha(t,\vec{x},\vec{k}(r,\mu,\varphi)).$$
We also obtain a transformed Boltzmann equation in divergence form for
our new PDF $\Phi$ in the $(x,y,z;r,\mu,\varphi)$ space, which reads
\begin{eqnarray}
\frac{\partial \Phi}{\partial t} &+&
\frac{\partial \mbox{ }}{\partial x} \left( a_{1} \, \Phi \right) +
\frac{\partial \mbox{ }}{\partial y} \left( a_{2} \, \Phi \right) +
\frac{\partial \mbox{ }}{\partial z} \left( a_{3} \, \Phi \right) +
\nonumber\\
&& \frac{\partial \mbox{}}{\partial r} \left( a_{4} \, \Phi \right)  + 
\frac{\partial \mbox{ }}{\partial \mu} 
\left( a_{5} \, \Phi \right) + \frac{\partial \mbox{ }}{\partial \varphi} \left(a_{6} \, \Phi \right) 
= C(\Phi), \nonumber
\end{eqnarray}
where the transport coefficients are, for
${(a_1,a_2,a_3) \propto \nabla_{\vec{k}} \varepsilon(\vec{k})}$, proportional to the
$\vec{k}$-gradient in transformed coordinates, and the rest are given by
\begin{eqnarray}
{a_{4}} &=&
{ - 2 \, \cd{E} \, \sqrt{r} \, \hat{e}_r\cdot \bEd } = - 2 \, \cd{E} \, \sqrt{r} \left(\mu,\sqrt{1 - \mu^{2}} \cos \varphi, \sqrt{1 - \mu^{2}} \sin \varphi\right)\cdot \bEd,\\
{a_{5} } &=& { - \cd{E} \frac{\sqrt{1 - \mu^{2}}}{\sqrt{r}} \, \hat{e}_{\mu}\cdot \bEd }
 =  - \cd{E} \frac{\sqrt{1 - \mu^{2}}}{\sqrt{r}} \left( \sqrt{1 - \mu^{2}}, -\mu \cos \varphi, - \mu \sin \varphi\right) \cdot \bEd ,\nonumber\\
{a_{6}} &=& {  - \cd{E} \, \frac{1}{\sqrt{r} \, \sqrt{1 - \mu^{2}}} \, \hat{e}_{\varphi} \cdot \bEd }
 =  - \cd{E} \, \frac{1}{\sqrt{r} \, \sqrt{1 - \mu^{2}}} \left(0, -\sin \varphi, \cos \varphi \right) \cdot \bEd .
\end{eqnarray} 

Regarding the transformed linear collision operator, we write $\vec{x} =
(x,y,z)$ and $\bfr = (r,\mu,\varphi)$, and obtain
\begin{eqnarray}
\mbox{} 
C(\Phi)(t,\vec{x},\bfr) =  
  \left. 
\frac{\sqrt{r}}{2} \, \int_{\Omega}  {\cal S}(\bfr', \bfr)
\right.
\Phi(t,\vec{x},\bfr') \: d\bfr'  
 - \, \Phi(t,\vec{x},\bfr) \int_{\Omega}
{\cal S}(\bfr, \bfr') \, \frac{\sqrt{r'}}{2} \,
  \: d\bfr' ,
  \nonumber
\end{eqnarray} 
showing the importance of the transformed PDF $\Phi$. Here ${\cal
  S}(\bfr', \bfr)$ represents the electron-phonon scattering for the
two-band system.







We use the dimensionless Poisson equation
\begin{equation}
 \nabla_{\vec{x}}\cdot \left( \epsilon_{r} \nabla_{\vec{x}} \nV
 \right)
 = 
\cd{p} \left[\rho(t,\vec{x}) - \mathcal{N}_{D}(\vec{x})\right],
\end{equation}
where 
\begin{equation} 
     \rho(t,\vec{x}) =
\int_{\Omega} \Phi_0(t,\vec{x},\bfr') \: d \bfr' ,
\qquad
\mathcal{N}_{D}(\vec{x}) =
    \left( \frac{\sqrt{2 \,\mass \kbt }}{\hbar} \right)^{\! \! -3}
    N_{D}( \vec{x}).
\end{equation}
The electron density is given by the first PDF coefficient, which
represents the mean of the PDF.

 
The discontinuous Galerkin method for the Boltzmann-Poisson system
represents a dynamic extension of the Gummel iteration map.  Starting
with an initial condition $\Phi_h$ and given boundary conditions, the DG
algorithm advances from $t^{n}$ to $t^{n+1}$ in these steps:
\begin{description}
\item[Step 1] Compute the charge density $\rho$.
\item[Step 2] Use $\rho$ to solve the Poisson equation (either by an
  integral form in 1D or by the LDG method in 2D or 3D) for the
  potential and electric field, and compute the transport coefficients
  $a_i$, $1\leq i\leq 6$.
\item[Step 3] Solve the transport part of the Boltzmann equation by
  DG, then use the method of lines for $\Phi_h$ (ODE system).
\item[Step 4] Evolve the ODE system by proper time stepping from
  $t^{n}$ to $t^{n+1}$ (if partial time steps are necessary, as in a
  Runge-Kutta method, repeat steps~1 to~3 as needed).
\end{description}

We use a rectangular Cartesian grid
in the transformed phase space. It has the form
$$
 \Omega_{ijkmn} = \left[ x_{i - \ot} , \, x_{i + \ot} \right] \times
                  \left[ y_{j - \ot} , \, y_{j + \ot} \right] \times
                  \left[ r_{k - \ot} , \, r_{k + \ot} \right] \times
                  \left[ \mu_{m - \ot} , \, \mu_{m + \ot} \right] \times
                  \left[ \varphi_{n - \ot} , \, \varphi_{n + \ot} \right]
$$
with
$1\leq i\leq N_x,\, 1\leq k\leq N_r,\, 1\leq m\leq N_{\mu},$ and
$ x_{i \pm \ot} = x_{i} \pm {\Delta x_{i}}/{2}$,
$ r_{k \pm \ot} = r_{k} \pm {\Delta r_{k}}/{2}$,
$ \mu_{m \pm \ot} = \mu_{m} \pm {\Delta \mu_{m}}/{2}$.

The test function $\psi(x,y,r,\mu,\varphi) \in {V_h}$ belongs to the
set of piecewise linear polynomials so that the set of all test
functions is
$$V_h := V_h^l := \left\{ v : v|{\Omega_{ijkmn}} \in
  P(\Omega_{ijkmn}^l) \right\},$$
where the $P(\Omega_{ijkmn}^l)$ are the
polynomials of degree $l\le 1$ on $ \Omega_{ijkmn}$.

 

Inside the cell $\mathring{\Omega}_{I}$, $I = (i,j,k,m,n)$, we
approximate the weighted PDF $\Phi$ by a linear polynomial in $V_h$,
i.e.,
\begin{equation}
\Phi_h   = 
  T_{I}(t) +
  X_{I}(t) \, \frac{(x - x_{i})}{\Delta x_{i}/2} +
  Y_{I}(t) \, \frac{(y - y_{j})}{\Delta y_{j}/2} +
  R_{I}(t) \, \frac{(r - r_{k})}{\Delta r_{k}/2} +
  M_{I}(t) \, \frac{(\mu - \mu_{m})}{\Delta \mu_{m}/2} +
  P_{I} \, \frac{(\varphi - \varphi_{n})}{\Delta \varphi_{n}/2}.
\nonumber 
\end{equation}
%
%
%
%
%
%
The charge density for a piecewise linear PDF~$\Phi$ is given by
\[
\rho(t,x,y) = 
\sum_{m=1}^{N_{\mu}} \sum_{n=1}^{N_{\varphi}}
\sum_{k=1}^{N_{r}} 
\left[ T_{ijkmn}(t) +
  X_{ijkmn}(t) \, \frac{(x - x_{i})}{\Delta x_{i}/2} +
  Y_{ijkmn}(t) \, \frac{(y - y_{j})}{\Delta y_{j}/2} 
\right] \Delta r_{k} \Delta \mu_{m} \, \Delta \varphi_{n}.
\]

In summary, the discontinuous Galerkin formulation for the vector
Boltzmann equation is to find $\Phi_h$ in the piecewise polynomial
space $V_h$ such that the equation
\begin{eqnarray} 
&& \int_K \frac{\partial \Phi_h}{\partial t} v_h d\Omega - \int_K \frac{\partial v_h}{\partial x} \left( a_{1} \, \Phi_h \right) d\Omega - \int_K \frac{\partial v_h}{\partial y} \left( a_{2} \, \Phi_h \right) d\Omega - \int_K \frac{\partial v_h}{\partial z} \left( a_{3} \, \Phi_h \right) v_h d\Omega  \nonumber \\
&& \mbox{ }  - \int_K \frac{\partial v_h}{\partial r} \left( a_{4} \, \Phi_h \right) d\Omega - \int_K \frac{\partial v_h}{\partial \mu} \left( a_{5} \, \Phi_h \right) d\Omega - \int_K \frac{\partial v_h}{\partial \varphi} \left(a_{6} \, \Phi_h \right) d\Omega  \nonumber \\ 
&& \mbox{} + F_{x}^{+} - F_{x}^{-} 
+ F_{y}^{+} - F_{y}^{-}  + F_{z}^{+} - F_{z}^{-} 
+ F_{r}^{+} - F_{r}^{-} + F_{\mu}^{+} - F_{\mu}^{-} + F_{\varphi}^{+} - F_{\varphi}^{-} = \int_K C(\Phi_h) v_h d\Omega, \nonumber 
\end{eqnarray}
where the $F^{\pm}$ represent the boundary integrals, holds for any
test function $v_h \in V_h $ and for each element $K =
\Omega_{ijkmn}$.


\section{SDG-BP: Stochastic Discontinuous Galerkin Method for the
  Boltzmann-Poisson System}

\subsection{The Symmetric Case: One-Dimensional in $x$, Two-Dimensional in $\vec{k}(r,\mu)$}

We consider a 1D $n^+$--$n$--$n^+$ silicon diode, rendering the
problem one-dimensional in position space.  The length of the diode is
$L=1 \,\mu\mathrm{m}$, and the length of the $n$-channel in the middle
is $400\,\mathrm{nm}$.  The doping concentration is $n^+ = 5\cdot
10^{23}/{\mathrm{m}^3} = 5\cdot 10^{17}/{\mathrm{cm}^3}$ and
$n=2\cdot10^{21}/{\mathrm{m}^3} = 2\cdot10^{15}/{\mathrm{cm}^3}$.

We consider the case with $k$-space azimuthal symmetry on $ \varphi
\in [0,2\pi]$ $\rightarrow \vec{k} = \vec{k}(r,\mu)$. Therefore, by the
symmetry assumptions, it is only necessary to consider the radial and
polar coordinates of the momentum.

The computational domain is taken as $ x \in [0,1]$, $r \in
[0,r_{max}]$, and $\mu \in [-1,1] $.  The constant $r_\mathrm{max}$ is
the cut-off such that $ \Phi(t,x,r,\mu) \approx 0$ for $ r \geq
r_\mathrm{max} $ in the numerical experiments. For example,
$r_\mathrm{max}\approx36$ for $V_\mathrm{bias} = 0.5\,\mathrm{V}$ in a
$400\,\mathrm{nm}$ channel.

The initial condition is $(\Phi_0,\Phi_1)(0,x,r,\mu) = (C N_D(x)
e^{-\varepsilon(r)} \sqrt{r}/2, 0)$ with a constant $C$ such that
$\rho(0,x) - N_D(x) = 0 $ at the initial time $t=0$.

The boundary conditions are the following.
\begin{itemize}
\item In the $x$-space, the charge concentration is neutral at the
  source and drain  endpoints $0 =x_{1/2}$ and $x_{N_x + 1/2}=1$.
This charge neutrality condition is imposed by  $\Phi(0, \vec{k}, t) = N_D (x) \frac{\Phi(x_1, \vec{k},
    t)}{\rho(x_1,\vec{k},t)}$ and
  $\Phi(1, \vec{k}, t) = N_D (x) \frac{\Phi(x_{N_x}, \vec{k}, t)}{\rho(x_{N_x},\vec{k},t)}$.
\item The applied potential (bias) is $V(0,\vec{k},t)=0$ and
  $V(1,\vec{k},t)=V_0$.
\item In the $(r,\mu)$-space, a cut-off boundary is used such that
  $\Phi$ vanishes at $r=r_\mathrm{max}$.
\item At the ``point'' boundaries, no boundary conditions are needed
  and transport equals zero analytically.  Hence, at the origin $r=0$,
  $a_4 = 0$ holds, and likewise at the poles $\mu = \pm 1$, $a_5 = 0$
  holds.  Therefore the boundary integrals are analytically equal to
  zero at $r=0$ and $\mu=\pm 1$.
\end{itemize}

Regarding time evolution, an RK2 Method was used in our simulations.

\subsection{Numerical Results}

\begin{figure} 
{\includegraphics[angle=0,width=1.0\linewidth]{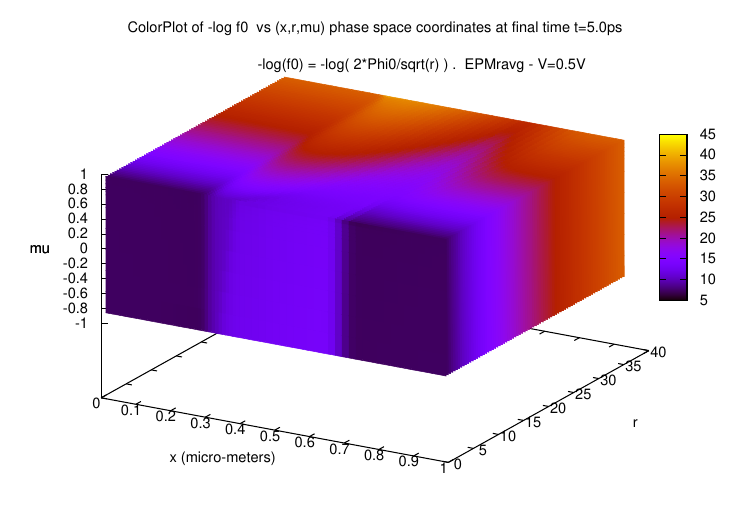}}
\caption{$-\log \alpha_0(x,r,\mu)$ for a
$1\,\mu\mathrm{m}$ diode, $0.5\,\mathrm{V}$ bias, and
$t_0=10.0\,\mathrm{ps}$.} 
\label{figs:NoRecomb0}
\end{figure}

\begin{figure} 
{\includegraphics[angle=0,width=1.0\linewidth]{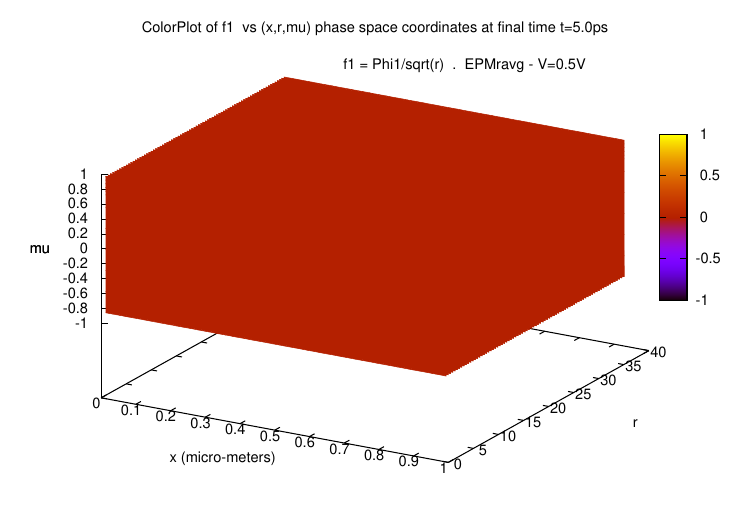}}
\caption{Coefficient $\alpha_1(x,r,\mu) = 0$ for a
$1\,\mu\mathrm{m}$ diode, $0.5\,\mathrm{V}$ bias, and
$t_0=10.0\,\mathrm{ps}$.} 
\label{figs:NoRecomb1}
\end{figure}

We present the numerical results for the coefficients of the truncated
PDF with a random variable.  We first do so for the benchmark case of
no recombination DG-BP: $-\log \alpha_0(x,r,\mu)$ (and
$\alpha_1(x,r,\mu)=0$ plotted directly) are shown for a
$1\,\mu\mathrm{m}$ diode, $0.5\,\mathrm{V}$ bias, and
$t_0=10.0\,\mathrm{ps}$. In this deterministic case, the reason why
there is no recombination lies in the vanishing of the first coefficient $\alpha_1 = 0$
related to random effects (Figure~\ref{figs:NoRecomb1}), and the zeroth
coefficient $\alpha_0$ contains all the information of the PDF
(Figure~\ref{figs:NoRecomb0}).

\begin{figure}
{\includegraphics[angle=0,width=1.0\linewidth]{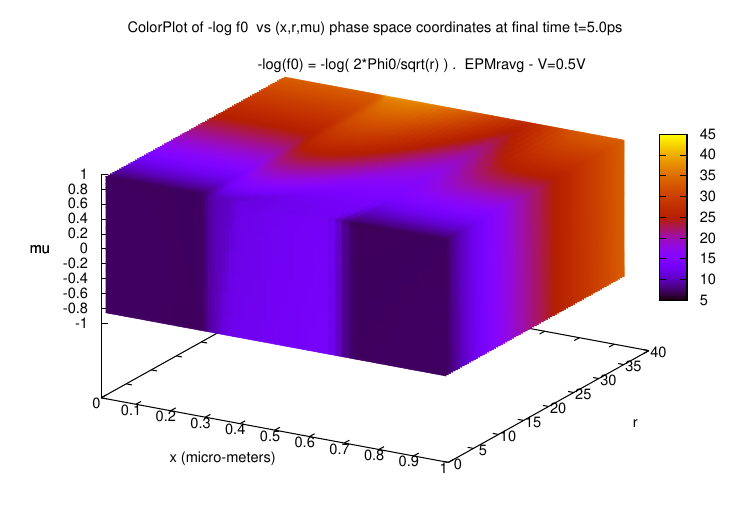}}
\caption{$-\log \alpha_0(x,r,\mu)$ for a
$1\,\mu\mathrm{m}$ diode, $0.5\,\mathrm{V}$ bias, and
$t_0=10.0\,\mathrm{ps}$.} 
\label{figs:Recomb0}
\end{figure}

\begin{figure}
{
\includegraphics[angle=0,width=1.0\linewidth]{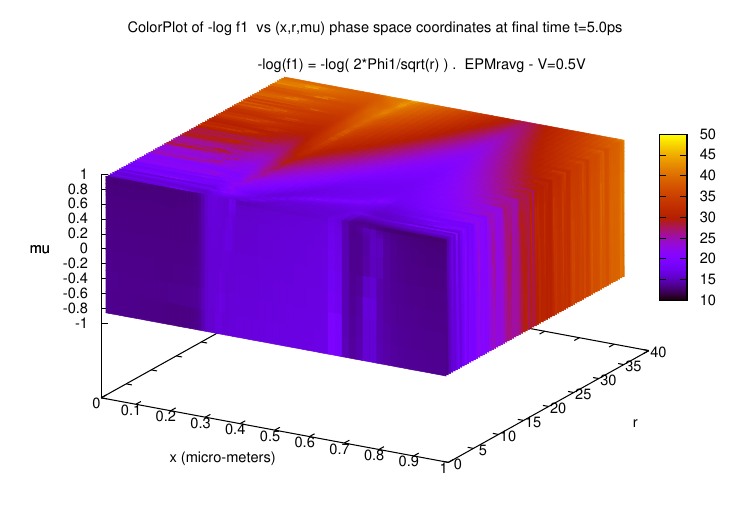} 
\caption{$-\log \alpha_1(x,r,\mu)$ for a
$1\,\mu\mathrm{m}$ diode, $0.5\,\mathrm{V}$ bias, and
$t_0=10.0\,\mathrm{ps}$.}
\label{figs:Recomb1} 
}
\end{figure}

We then consider the PDF coefficients from the simulations of the
SDG-BP system with the recombination terms $-\log \alpha_0(x,r,\mu,t)$
in Figure~\ref{figs:Recomb0} and $-\log \alpha_1(x,r,\mu,t)$ in
Figure~\ref{figs:Recomb1} for a $1\,\mu\mathrm{m}$ diode,
$0.5\,\mathrm{V}$ bias, and $t_0=10.0\,\mathrm{ps}$ as well. The
variations in $\alpha_1$ are located in similar regions of the phase
space, while for $\alpha_0$ they seem finer and more pronounced.

\begin{figure}
{\includegraphics[angle=0,width=1.0\linewidth]{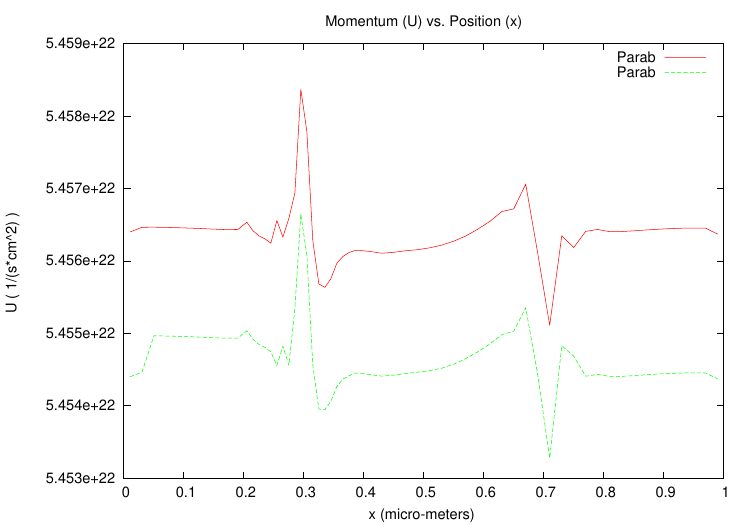}}
\caption{Momentum (current) $M(x,t)$.}
\label{fig:CompareMomentum}
 
\end{figure}

We also compare SDG-BP with recombination terms against SDG-BP with no
recombination case by calculating the moments with $\Phi_0$ for both.  The difference is observed mainly in the prediction of the
momentum (current) two orders of magnitude below the mean value of the
current (Figure~\ref{fig:CompareMomentum}), indicating the finer
resolution of the momentum by use of the stochastic Galerkin method.
We also plot the expectation, variance, and standard deviation of our
probability density function in the SDG-BP method, given as
$\mathrm{E}[f] = \alpha_1$, $\mathrm{Var}[f]=\sum_{k=2}^2 \alpha_k^2 =
\alpha_2^2$, and $\mathrm{S}[f]= \sqrt{\sum_{k=2}^2 \alpha_k^2} = {
  |\alpha_2| }$, respectively.

\begin{figure}
\includegraphics[angle=0,width=1.0\linewidth]{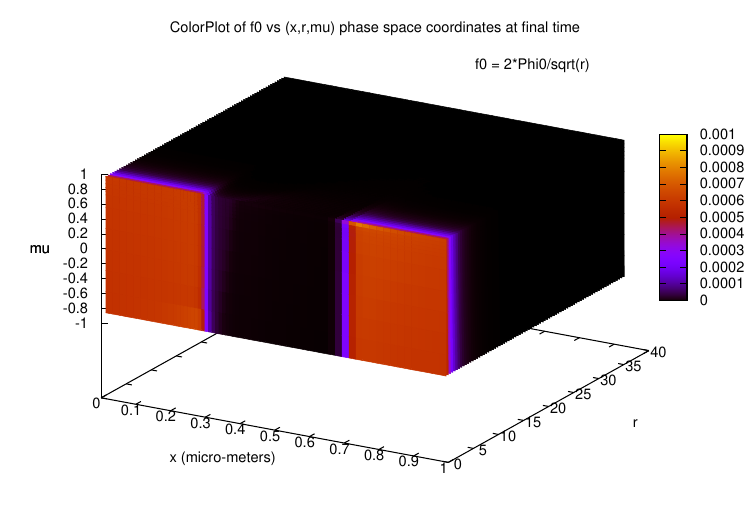} 
\caption{Coefficient $\alpha_0(x,r,\mu)$ for a
$1\,\mu\mathrm{m}$ diode, $0.5\,\mathrm{V}$ bias, and
$t_0=10.0\,\mathrm{ps}$.}
\end{figure}

\begin{figure}
\includegraphics[angle=0,width=1.0\linewidth]{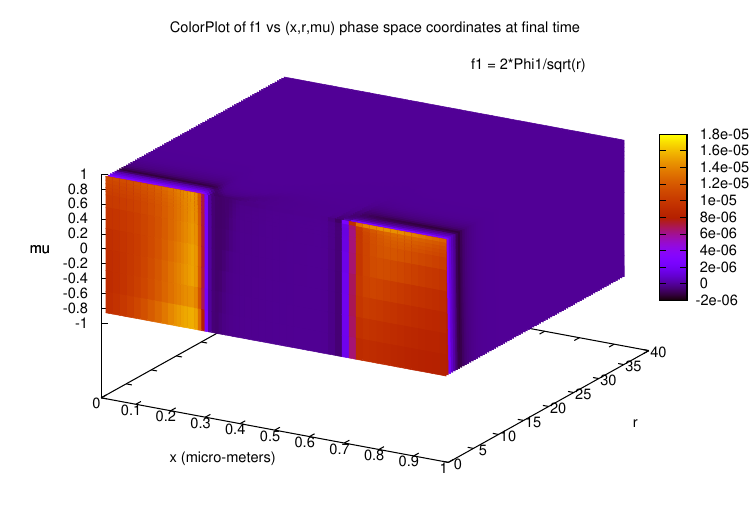} 
$\quad$
\caption{Coefficient $\alpha_1(x,r,\mu)$ for a
$1\,\mu\mathrm{m}$ diode, $0.5\,\mathrm{V}$ bias, and
$t_0=10.0\,\mathrm{ps}$.}
\end{figure}

\begin{figure}
\includegraphics[angle=0,width=1.0\linewidth]{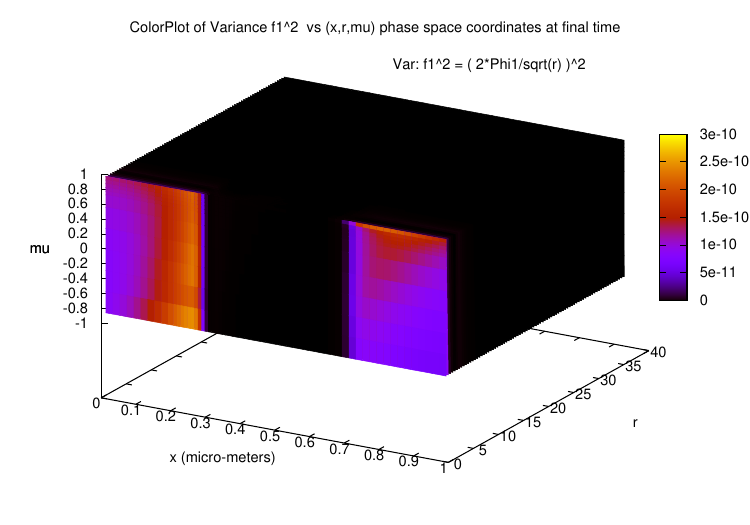} 
\caption{Variance $\mathrm{Var}[f]= \alpha_2^2$.}
\end{figure}

\begin{figure}
\includegraphics[angle=0,width=1.0\linewidth]{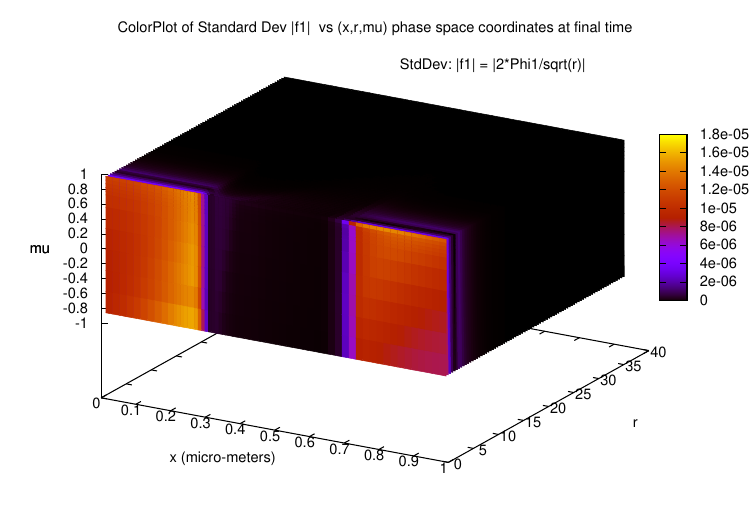} 
$\quad$
\caption{Standard deviation $\mathrm{S}[f]= |\alpha_2|$.}
\end{figure}

\begin{figure}
{
\includegraphics[angle=0,width=1.0\linewidth]{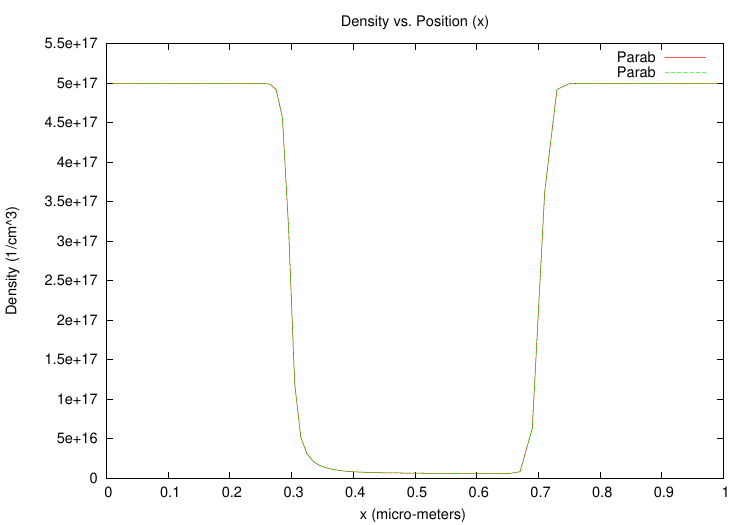} 
\caption{Density $\rho(x,t)$.}
} 
\label{fig:CompareDensity}
\end{figure}

\begin{figure}
{
\includegraphics[angle=0,width=1.0\linewidth]{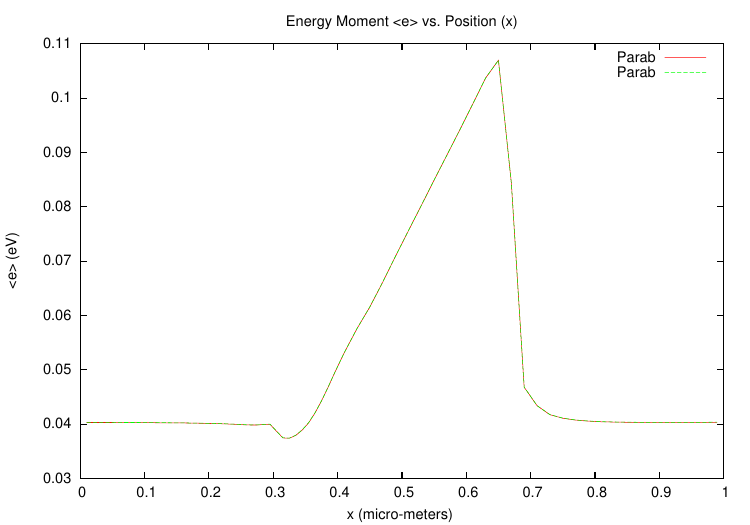} 
\caption{Energy $e(x,t)$ for a
$1\,\mu\mathrm{m}$ diode, $0.5\,\mathrm{V}$ bias, and
$t_0=10.0\,\mathrm{ps}$.} 
\label{fig:CompareEnergy}
}
\end{figure}

\begin{figure}
{
\includegraphics[angle=0,width=1.0\linewidth]{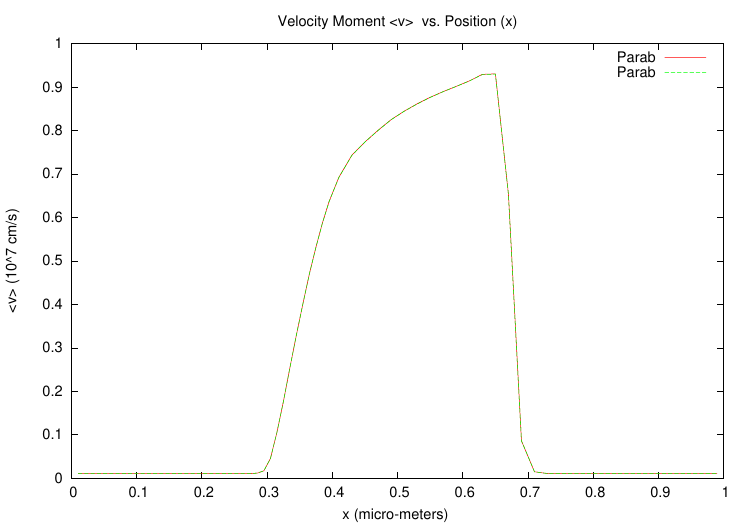} 
\caption{Velocity $v(x,t)$.}
\label{fig:CompareVelocity} 
}
\end{figure}

\begin{figure}
{\includegraphics[angle=0,width=1.0\linewidth]{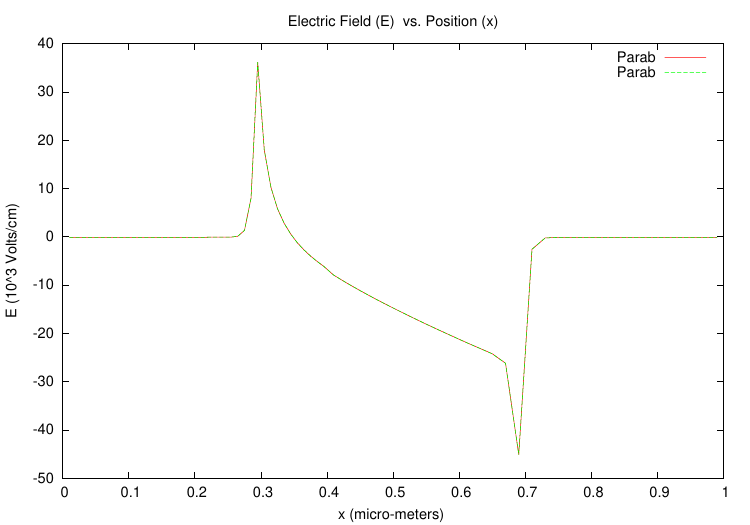}}
\caption{Electric field $\mathrm{E}(x,t)$.} 
\label{fig:CompareEfield}
\end{figure}

\begin{figure}
{\includegraphics[angle=0,width=1.0\linewidth]{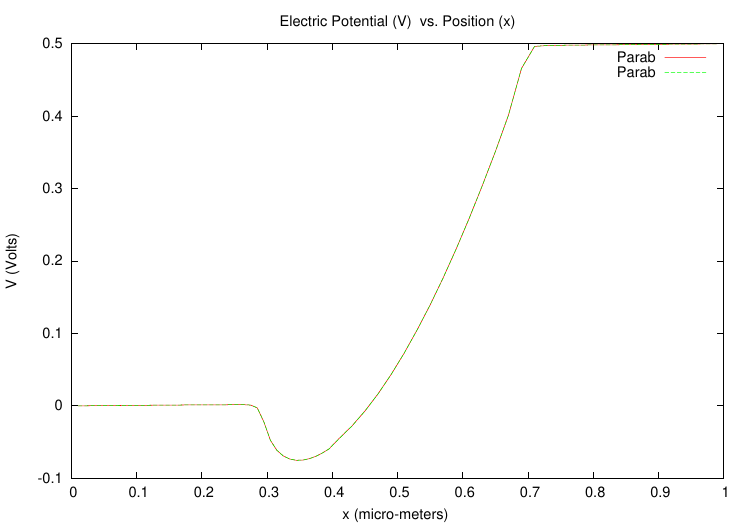}}
\caption{Potential $V(x,t)$.} 
\label{fig:ComparePotential}
\end{figure}

\section{Conclusions}

Uncertainty quantification in the Boltzmann-Poisson system is crucial
by the own probabilistic nature of the problem due to the high number
of particles involved and due to its quantum mechanical features.
Studying randomness in the temperature is an important leading
example, both for reasons stemming from the physical nature of the
problem -- as the environment temperature may fluctuate -- as
well as related to the mathematical aspects, since the temperature is
a scalar random variable that introduces randomness in the collision
term 
and more precisely in the coefficients 
multiplying the Dirac delta
distributions 
appearing in the electron-phonon collisions by the Fermi golden rule.

Our numerical results for the stochastic Galerkin method for the
Boltzmann-Poisson system, assuming a random temperature in an
electron-phonon collision operator, show a coefficient $\alpha_1$
related to randomness, whose variations are located in similar regions
of the phase space to the ones of the average term $\alpha_0$, which has no
randomness but shows finer and more pronounced variations.

Our comparison of SDG-BP with recombination terms against the
no-recombination case, using the mean term $\alpha_0$ in both cases to
calculate the moments, shows only a difference in the moment (current)
between these two cases. We have to remember that the moment is the
product of density and energy, which on their own scales do not seem
to exhibit a large difference between these two cases. However, the
moment has a value close to a constant over the position domain when
equilibrium is reached, and the difference in this mean value is
observed between the two cases, although it is two orders of magnitude
below the average value of the current
(Figure~\ref{fig:CompareMomentum}).  Here a truncated random expansion
up to first order in~$z$ was employed, which therefore discards terms
of order $z^2$ whose average is non-zero over the domain.  This is the
likely explanation for the slight difference in the momentum values
between the random and deterministic simulations.

We have also devised numerical methods for quantifying
uncertainty related to the phonon energy via stochastic Galerkin
methods, which can be handled by the introduction of
distributional derivatives with respect to the random variable. This
mathematical structure departs from the usual form of the collision
term in stochastic Galerkin for Boltzmann models by the need of
distributional derivatives in the random space, being the first case
in stochastic Galerkin methods for kinetic equations where this
structure appears, and opening a new analytical treatment of randomness in the aforementioned stochastic method. 

In conclusion, 
we have handled in this work the possible uncertainties arising in a model of electron transport in semiconductors by stochastic Galerkin, mainly related to the collision mechanisms in this paper.
We calculate the propagated uncertainty in the electron probability density function due to possible uncertainties in either the phonon energy (adding a random variable given by either a Gaussian or uniform distribution, considering first an approximate randomness to first order in the phonon energy and then the full calculation) or in the lattice temperature (assumed to vary randomly according to a uniform distribution). 

Our purpose is to observe how physical variables which can either
behave randomly in a real world setting (such as a varying
temperature) or are known to be described approximately in our
model (such as the phonon energy, which is often assumed to be
constant, but really is known experimentally to be not constant) can
affect physical observables such as electric current, average energy
or density, since our kinetic model lets us calculate those measurable
quantities by means of moments of the PDF.

Our study is useful to let us predict in real world settings the
impact that uncertainties or limitations in commonly used idealized
models have on the behavior of an electronic device such as a diode or
a MOSFET. This study is also useful in terms of introducing
uncertainties in the energy transition arguments of the collision
integrals with the ultimate goal of jumping from the scalar treatment
presented in this paper to the case of an energy band structure
$\varepsilon(\vec{k})$, which is a scalar function of a vector
variable, in future work.

We have calculated with our numerical methods the variation in kinetic
moments (density, mean energy, current, etc.)\ associated with a
physically reasonable temperature variation in the lattice
environment. We have presented the algorithms to make an approximate
(to first order) and exact calculation of the propagation of
uncertainty of the phonon energy (bounding the error of a constant
phonon energy by a uniform distribution) into the PDF which will give
as well the associated uncertainty in the prediction of moments.



\section{Acknowledgment}

This work was supported by FWF (Austrian Science Fund) START project
Y660 \textit{PDE Models for Nanotechnology}.  The first author also
acknowledges start-up support from The University of Texas at San
Antonio.

\bibliography{myreferences} 
\bibliographystyle{ieeetr}

\end{document}